\documentclass{elsart}
\usepackage{natbib}

\usepackage[latin1]{inputenc}
\usepackage{amssymb, amsmath}
\usepackage{epsfig, epic, eepic}

\def\GL{\mathrm{GL}}

\def\J{\mathrm{J}}
\def\x{\mathrm{x}}
\def\a{\mathrm{a}}

\def\S{\mathcal{S}}

\def\Bar{\overline}

\def\fbar{{\bar{f}}}
\def\xbar{{\bar \x}}
\def\E{\mathcal{E}}

\def\S{\mathcal{S}}

\newcommand{\Frac}[2]{\displaystyle \frac{#1}{#2}}
\newcommand{\Sum}[2]{\displaystyle{\sum_{#1}^{#2}}}

\newcommand{\N}{{\mathbb N}}

\newcommand{\C}{{\mathbb C}}

\def\d{\mathrm{d}}
\def\J{\mathrm{J}}
\def\I{\mathrm{I}}
\def\x{\mathrm{x}}
\def\a{\mathrm{a}}

\def\mytilde{\hspace{-0.1em}\raisebox{-.8ex}{\textasciitilde}\hspace{-0.1em}}

\newtheorem{definition}{Definition}
\newtheorem{theorem}{Theorem}
\newtheorem{proposition}{Proposition}

\newtheorem{proof}{Proof}
\newtheorem{example}{Example}

\newcommand{\pref}[1]{(\ref{#1})}

\newenvironment{system}[1][]%
	{\begin{eqnarray} #1 \left\{ \begin{array}{lll}}%
	{\end{array} \right. \end{eqnarray}}

\begin{document}

\begin{frontmatter}

\title{Elie Cartan's Geometrical Vision or How to Avoid
	Expression Swell}

\author{Sylvain Neut},
\ead{neut@lifl.fr}
\author{Michel Petitot},
\ead{petitot@lifl.fr}
\author{Raouf Dridi}
\ead{dridi@lifl.fr}
\address{Universit\'e Lille I,
	LIFL, bat. M3,
	59655 Villeneuve d'Ascq CEDEX,
	France}

\begin{abstract}
The aim of the paper is to demonstrate the superiority of Cartan's
method over direct methods based on differential elimination for
handling otherwise intractable equivalence problems. In this sens,
using our implementation of Cartan's method, we establish two new
equivalence results. We
establish when a system of second order ODE's is equivalent to flat
system (second derivations are zero), and when a system of holomorphic
PDE's with two independent variables and one dependent variables is
flat. We consider the problem of finding transformation that brings a
given equation to the target one.  We shall see that this problem becomes
algebraic when the symmetry pseudogroup of the target equation is zerodimensional. We avoid the
swelling of the expressions, by using  non-commutative derivations
adapted to the problem.

\end{abstract}

\begin{keyword}
 Cartan's equivalence method  \sep  differential algebra
\end{keyword}

\end{frontmatter}

\section*{Introduction}
Present  ODE-solvers  make use of a 
combination of symmetry methods and classification methods. 
Classification methods are used when the ODE matches
a recognizable pattern 
(e.g. as listed in \citep{kamke}).
Significant progress would be made if it was possible
to compute in advance the differential invariants  
that allow to decide whether the  equation to be  solved is equivalent to one of
the list by a change of coordinates. We will show that,
for the computation of these invariants, the  geometrical approach 
offers advantages over non geometrical approaches (e.g.  Riquier, Ritt, Kolchin etc.)

The  change of coordinates that maps a given differential equation to a target one is a  solution of certain PDE's system. 
Differential algebra allows us to compute the integrability conditions of this  system.
Unfortunately, in practice one is often confronted with computer output consisting of
several pages of  intricate formulae. Even so, for more complicated
 problems with higher complexity such those treated in
section 5, and  many others coming from biology, physics,
etc. differential  algebra is not efficient due to  expression swell.

 In his equivalence method, \'Elie Cartan  formulated the PDE's system as
 a linear Pfaffian system. In this way, the integrability conditions
 are computed with a process  called absorption of torsion leading to
 sparse structure equations.  In addition, this computation is done by
 separately and  symmetrically treating the linear Pfaffian
 system. This divides  the number of variables by two.

Regarding ODE-solvers, Cartan's equivalence method is complementary to  symmetry
method. In the case when the  symmetry pseudogroup of 
 the input equation  is
zerodimensional, actual DE-solvers are unusable.
Thus, in such a case,  one can map this equation to a known equation. 
We shall see that the  change of coordinates realizing the equivalence
can be computed without integrating differential equation. 

The paper is organized as follows. In section 1, we define the
equivalence problem of differential equations under a given
pseudogroup of diffeomorphisms. In section 2, using {\sc Diffalg}, we
solve two examples which will serve  for our
argumentation and we finger out the limitations of such technique. We
give a first optimization, to avoid the expression swell, by using 
non commutative derivations. The goal of section 3 is to introduce the
reader to the calculation of the integrability conditions of linear Pfaffian
 systems. In section 4, Cartan's method is applied to the
equivalence  problem  presented in section 1. In section~5, we give
new results on the equivalence of ODE's  and PDE's systems with flat
systems (Theorems \ref{SysEdoThm1} and \ref{thm1}). In section 6,  
we consider the problem of finding transformation that brings a
given equation to the target one which is usually even harder than, to
establish the equivalence. We shall see that this problem becomes
algebraic when the symmetry pseudogroup is zerodimensional. We avoid the
swelling of the expressions, by using 
 non-commutative derivation adapted to the problem. 


\section{Formulation of local equivalence problems}
An  equivalence problem is the following data : a class of (systems of)
differential equations~$\E_f$ and a pseudogroup of transformations
acting on this class.  The two  differential equations $\E_f$ and $\E_{\bar f}$ are
 said to be equivalent  under the pseudogroup~$\Phi$
 (we write $\E_f \sim_{\Phi } \E_{\bar f}$)   if and only if there
 exists a  local diffeomorphism $\varphi \in \Phi$ which maps the
 solutions  of~$\E_f$ to the solutions of $\E_{\bar f}$.
The change of coordinates $\varphi \in \Phi$ 
is solution of a PDE's system that we can generated by an algorithm.
A program like {\sc Diffalg} \citep{BLOP, boulier} or {\sc Rif}
\citep{Reid96}  can  compute the integrability conditions
of such  PDE's system, then the existence of $\varphi$ is decidable.

Consider the two second order ODE $(\E_f)$ and $( \E_{\bar f})$
\begin{eqnarray}  \label{exemple}
	\Frac{\d^2 y}{\d x^2} = f\left(x,\, y,\, \Frac{\d y}{\d x}\right)
	&\mbox{ and }&
	\Frac{d^2 \bar y}{\d \bar x^2} =
	\bar f\left(\bar x,\, \bar y,\, \Frac{\d \bar y}{\d \bar x}\right).
\end{eqnarray}
Let  $\Phi$ denote  the pseudogroup of local diffeomorphisms
$\varphi: \C^2 \to \C^2$   defined by the Lie equations 
$$\Frac{\partial\bar x}{\partial x}=1, \quad  \Frac{\partial \bar x}{\partial y}=0,\quad \Frac{\partial \bar y}{\partial y} \neq 0. $$
This gives
$
	(\bar x,\, \bar y) = \varphi(x,y) = (x+C,\, \eta(x,y))
$
where $C$ is a constant and $\eta(x,y)$ is an arbitrary function.


\section{Equivalence problems and differential elimination}
We shall consider  equivalence problems
with fixed (determined) \emph{target} equation~$\E_{\bar f}$. The question is to find the
explicit  conditions on $f$ such that $\E_f \sim_{\Phi } \E_{\bar f}$.   

For instance, consider the equivalence problem of $2^{nd}$ order ODE
presented in the previous section. Let $(x,\,y,\,p=y',\,q=y'')$ denote
a local coordinates system of  the jet space
$\J^2(\C,\C)$. Thus,   the problem reads
\begin{equation}\label{existe}
        \exists\, \varphi \in \Phi, \quad \varphi^{*}(\bar{q}-\bar
        f(\bar x,\bar y,\bar p)) = 0 \ \mod \ q - f(x,y,p)=0.
\end{equation}
The prolongation formulae \citep{Olver} of $\varphi$ are
\begin{eqnarray}\nonumber
\begin{gathered}
 \varphi^{*}\bar{p} =D_x \eta = \eta_{x} +\eta_{y} p, \quad
 \varphi^{*}\bar{q} = D_x^2 \eta =  \eta_{xx} + 2\eta_{xy}p + \eta_{yy}p^2 + \eta_{y}q
\end{gathered}
\end{eqnarray}
where  $D_x = \frac{\partial}{\partial x}
             + p\ \frac{\partial}{\partial y}
             + f(x,y,p)\frac{\partial}{\partial p}$ is the total
             derivative.  
The following examples explain how one can use differential
elimination to solve such question.

\begin{example}
Let us suppose that the \emph{target} $\bar f$
is identically zero.  
The equations~\pref{existe} take the form of a polynomial PDE's system
\begin{equation}\label{syst_1}
        \eta_{xx} +  2\eta_{xy}p +  \eta_{yy}p^2 + \eta_{y}f = 0, \quad
        \eta_{p} = 0,  \quad
        \eta_y \neq 0.
\end{equation}
By eliminating $\eta$ in \pref{syst_1} using the ranking $\eta \succ f$,
we obtain the characteristic~set
\begin{equation} \label{etaf}
\begin{gathered}
   \eta_{xx}  =  -\eta_{y} f_{} + p \eta_{y} f_{p} - \Frac{1}{2} p  \eta_{y} f_{pp}, \quad
   \eta_{xy}  =  - \Frac{1}{2} \eta_{y} f_{p} + \Frac{1}{2} p \eta_{y} f_{pp},  \\
   \eta_{yy}  =  - \Frac{1}{2} \eta_{y} f_{pp}, \quad
   \eta_{p}   =  0, \\
   f_{ppp}    =  0, \quad
   f_{xp}     =  -f_{pp} f_{} + 2 f_{y} + \Frac{1}{2} f_{p}^2 - p f_{yp}.
  \end{gathered}
\end{equation}
It follows that 
the $2^{nd}$ order  differential equation  $y''=f(x,y,y')$ is reduced to the equation $\bar y''=0$
by a transformation~$\varphi$ of the form  
 $\varphi(x,y) = (x+C,\,\eta(x,y))$ 
if and only if
\begin{equation}\label{thm1}
\begin{gathered}
   f_{ppp}    =  0, \quad
   f_{xp}     =  -f_{pp} f_{} + 2 f_{y} + \Frac{1}{2} f_{p}^2 - p f_{yp}.
  \end{gathered}
\end{equation}
To obtain the change of coordinates $\varphi$ we have to integrate the PDE's
system given by the first four equations of \pref{etaf}. 
\end{example}

The next example shows that, in  favorable cases, we can
determine the change of coordinates without any integration.

\begin{example}
Suppose that the \textit{target} equation is  $(P_I)$, the   Painlev\'{e} first
equation 
$
\bar y''=6\bar y^2+\bar x
$
. The problem formulation is as above and in this case  {\sc Diffalg} returns the following  characteristic set
\begin{eqnarray}\label{sysp1}
\eta(x,y) &=& 1/12\,f_{{y}}-1/24\,f_{{x,p}}-1/24\,f_{{p,p}}f_{{}}+1/48\,{f_{{p}}}^{2}-1/24\,pf_{{y,p}},\\ \nonumber
C &= & -x+1/16\,pf_{{p}}f_{{y,p,p}}f_{{}}+1/16\,{p}^{2}f_{{p,p}}f_{{p}}f_{{y,p}}+1/16\,pf_{{p}}f_{{p,p}}f_{{y}}\\ \nonumber
~~&& \vdots \\ \nonumber
f_{{x,x,x,x,p}} & = & -24+5/2\,f_{{}}pf_{{p,p}}f_{{p}}f_{{y,p}}-4\,f_{{x}}f_{{x,y,p}}+2\,f_{{x,x,x,y}}-2\,pf_{{y,p}}f_{{p}}f_{{x,p}}\\ \nonumber
~~&& \vdots \\ \nonumber
f_{{x,x,y,p}} &=& -pf_{{p,p}}f_{{y,y}}+3\,pf_{{y,p,p}}f_{{y}}+f_{{p,p}}{p}^{2}f_{{y,y,p}}-3/2\,pf_{{p}}f_{{y,y,p}}\\ \nonumber
~~&& \vdots\\ \nonumber
f_{{x,y,y,p}} & = & 2\,f_{{y,y,y}}+{f_{{y,p}}}^{2}-pf_{{y,y,y,p}}-f_{{y,y,p,p}}f_{{}}-2\,f_{{y,p,p}}f{{y}}-\\ \nonumber
~~&& 1/2\,f_{{p,p}}f_{{y,p,p}}f_{{}}+f_{{p}}f_{{y,y,p}}-1/2\,f_{{p,p}}f_{{x,y,p}}-1/2\,f_{{p,p}}pf_{{y,y,p}}-\\ \nonumber
~~&& 1/2\,{f_{{p,p}}}^{2}f_{{y}}+1/2\,f_{{p,p}}f_{{p}}f_{{y,p}},\\ \nonumber
f_{{x,p,p}} & = & f_{{y,p}}-pf_{{y,p,p}},\\ \nonumber
f_{{p,p,p}} &=& 0. \nonumber
\end{eqnarray}%
The two first equations demonstrate that the change of variable $\varphi$ is obtained
without integrating differential equation and is unique. We shall  see,
in section~6,  that this results from
the fact that the  symmetries  pseudogroup of  $(P_I)$ is zerodimensional
(in fact reduced to the identity).
The other equations gives the requested conditions on $f$. 
\end{example}

As the reader may have noticed, such explicit formulae consisting of
several lines (pages) prove quite useless for practical
application. In section~6, we shall wee that the same formulae take a more compact form
when they are  written in term of the associated  invariants.

 More dramatically, the above brute-force  method is rarely efficient due to
expression swell involved by the use of commutative derivations. 
Indeed, the study of the
equivalence of the $3^{rd}$ order differential equation $y'''= f(x, y, y',y'' )$ with $\bar y'''=0$ under
contact transformations
$
 (\bar{x},\bar{y})  = \left (\xi(x,y,y'), \eta(x,y,y')\right),
$
requires the prolongation to $\J^3 = (x,\ y,\ p = y',\ q = y'',\
r=y''')$, that is to~find 
\begin{eqnarray*}
\bar p  &=&  \frac{D_{x}\eta}{D_{x}\xi}, \\
\bar{q} &=&  \frac{{D_{x}}^{2}\eta D_{x}\xi - D_{x}\eta{D_{x}}^{2}\xi}{{D_{x}\xi}^{3}},\\
\bar{r} &=&  \frac{ {D_{x}}^{3}\eta {D_{x}\xi}^{2}
                - 3 {D_{x}}^{2}\eta {D_{x}}^{2}\xi D_{x}\xi
                + 3 D_{x}\eta{{D_{x}}^{2}\xi}^{2} -  D_{x}\eta {D_{x}}^{3}\xi D_{x}\xi}{D_{x}\xi^5},
\end{eqnarray*}
where $D_x =\frac{\partial}{\partial x} + p\frac{\partial}{\partial y}
+  q \frac{\partial}{\partial p} + f(x,y,p,q) \frac{\partial}{\partial q}$.
The change of coordinates satisfies the PDE's system $\{\bar r=0, \ 
\frac{\partial \bar{p}}{\partial q } = 0,\  \xi_q = 0,\ \eta_q = 0\}$. 
Using commutative derivations, this system blows up and takes more than one hundreds of lines
\begin{eqnarray}\ \nonumber
\bar{r} & = &  3 p \eta_{xxy} {\xi_{x}}^{3} - 12 \eta_{xx} {p}^{2} \xi_{xy} \xi_{y} \xi_{x} - 6 {p}^{3} \eta_{y} \xi_{xxy} \xi_{y} \xi_{x} + 3 p \eta_{y} {\xi_{xx}}^{2} \xi_{x} +\\ \nonumber
&& \vdots \\  \nonumber
&& \textrm{\bf {100 lines of differential polynomials}} \\  \nonumber
\end{eqnarray}%
and calculation (treatment by {\sc Diffalg}) do not finish !

The first optimization is to use  the non
 commutative derivations $\{D_x, \frac{\partial}{\partial
   y},\frac{\partial}{\partial p}, \frac{\partial}{\partial q}\}$ such
 in \citep{hubert00} inspired by \citep{neut:these}. A better alternative is to use the associated invariant derivations discussed later on.

\section{The geometrical approach of  the
  integrability conditions calculation}
 \'E. Cartan's transforms an analytic PDE's system into an equivalent
\emph{linear} Pfaffian system (with a condition that specifies the independent variables).
He gave a method to compute the integrability conditions of any
analytic linear Pfaffian systems and therefore
of any  analytic PDE's systems. This algorithm is based on the process of
absorption of torsion which leads   to sparse equations.

Recall that a  Pfaffian system (with independence condition) on real
analytic  manifold $M$ is a an exterior differential system of the form
\begin{system}
\omega^\alpha =0, &&\quad (1 \leq \alpha \leq a)  \\   
\theta^1 \wedge \theta^2 \wedge \cdots \wedge \theta^n &\neq& 0,
\end{system}%
 where    $\omega^\alpha$ and  $\theta^i$ are linearly independent
 differential 1-forms defined on~$M$ and  $a, n\in \N$.
An \emph{integral manifold} $i: S \hookrightarrow M$  is a submanifold
 $S$ of $M$ such that $i^* \omega^\alpha =0$ for all
$1 \leq \alpha \leq a$ and
$i^* (\theta^1 \wedge \theta^2 \wedge \cdots \wedge \theta^n) \neq 0$.

 Let  $[I]$ and  $[J]$ denote the exterior differential ideals
 respectively generated  by $I=(\omega^\alpha)$ and $J=(\omega^\alpha,\, \theta^i)$ for
$1 \leq \alpha \leq a$ and $1\leq i \leq n$.
\begin{definition} \label{pfaff:lineaire}
    A Pfaffian system $I \subset J \subset \Omega^1 M$ is 
    \emph{linear} if and only if $\d I = 0 \mod [J]$.
\end{definition}

 One obtains a local basis $(\omega^\alpha, \, \theta^i, \, \pi^\rho)$
 of $\Omega^1 M$ by  completing the basis $(\omega^\alpha, \,
 \theta^i)$ of $J$ by  the 1-forms $\pi^\rho\in \Omega^1 M$ ($1\leq
 \rho \leq r$).  If the Pfaffian system is linear then there exist analytic
 functions $A^\alpha_{\rho i}$ and $T^\alpha_{jk}$ defined on $M$ such that
\begin{equation} \label{eq:structure}
     \d \omega^\alpha = A^\alpha_{\rho i}\ \pi^\rho \wedge \theta^i
     	+ \Frac{1}{2} T^\alpha_{jk}\ \theta^j \wedge \theta^k \mod [I],
     	\qquad (1\leq \alpha \leq a).
\end{equation}

\begin{proposition}
Given two analytic manifolds $X$ and $U$. Every $q$-order PDE's system $\E^q \subset \J^q(X,U)$ is equivalent to a
linear  Pfaffian systems defined on $\E^q$. 
\end{proposition} 
\begin{proof}
 Suppose that $x = (x^i)_{1\leq i \leq n:= \dim X}$ are local  coordinates of
 $X$ and  $u=(u^\alpha)_{1\leq \alpha \leq m := \dim U}$  are local
 coordinates of $U$, then $(x^i, u^\alpha_{I})$, where $I=i_1 i_2
 \cdots i_\ell$  and $0\leq l\leq q$, constitutes a local coordinates of
 $\J^q$.  The contact system of $\J^q$ is
 $\{\d u^{\alpha}_{I} - u^{\alpha}_{I,i} \d x^i = 0,\quad
		\ell < q,\quad  1\leq \alpha \leq m \}
$ 
with $\d x^1  \wedge \cdots \wedge \d x^n \neq 0$
and these equations continue to hold when we restrict to~$\E^q$.
\end{proof}

\subsection{Essential elements of torsion}
 Since 
$ i^*(\theta^1)\wedge i^*(\theta^2)\wedge \cdots \wedge i^*(\theta^n) \neq 0, $
the forms  $i^*(\pi^\rho)$ are linear combinations of the forms
$i^*(\theta^i)$, i.e.  there are coefficients $\lambda^\rho_i $ such
$ i^*(\pi^\rho) =  \lambda^\rho_i i^*(\theta^i)$. Substituting into
\pref{eq:structure} leads to 
 $$ \Sum{1\leq j<k \leq n}{ } i^ * \Big(T^\alpha_{jk } -  
A^{\alpha}_{\rho j } \lambda^\rho_k +
 A^{\alpha}_{\rho k } \lambda^\rho_j \Big) \ i^ * (\theta^j \wedge \theta^k) = 0, \qquad (1 \leq \alpha \leq a).  $$
These conditions are equivalent to the system (by omitting the pullback $i^*$)
\begin{equation} \label{eq:lambda}
	T^\alpha_{jk} = A^{\alpha}_{\rho j} \lambda^\rho_k
		- A^{\alpha}_{\rho k} \lambda^\rho_j,
		\qquad  (1\leq \alpha \leq a; 1\leq j < k \leq n)
\end{equation}
which is \emph{linear} in the unknown coefficients
$\lambda^\rho_i$. By eliminating the  coefficients~$\lambda^\rho_i$
using the standard  Gaussian elimination, one obtains linear
combinations of  the functions $T^\alpha_{jk}$, called
\emph{essential}  torsion elements, which  inevitably vanish.  
 Reciprocally, this vanishing ensures the existence of the $\lambda^\rho_i$.
\begin{theorem}  \label{cond:compatibilite}
   The essential torsion elements are real--valued functions defined
   on $M$ and vanishing  on any integral manifold   $S \subset M$ of
   the linear  Pfaffian system $I \subset J \subset \Omega^1 M$. In
   other words,  they are the integrability conditions.
\end{theorem}
\subsection{Absorption of torsion}
Elie Cartan was used to calculate the integrability conditions
provided by the  preceding theorem using a process,  called today \emph{absorption of torsion}.   
In the structure equations  \pref{eq:structure}, let us replace the
$\pi^\rho$ by the general linear combination
$	\Bar \pi^\rho := \pi^\rho + \lambda^\rho_i\ \theta^i.$ 
 This yields 
\begin{equation}  \label{formules:chgt:base}
    \Bar A^{\alpha}_{\rho i} = A^\alpha_{\rho i},  \quad
    \Bar T^{\alpha}_{jk} = T^{\alpha}_{jk}
    		+ A^{\alpha}_{\rho j} \lambda^\rho_k
    		- A^{\alpha}_{\rho k} \lambda^\rho_j.
\end{equation}%
The process of absorption of torsion consists in calculating 
the $\lambda^\rho_i$   so that to fix the maximum of  $\bar T^{\alpha}_{jk}$ to zero.
After the  absorption, the torsion elements  remaining  non zero form
a  basis of the essential torsion elements.
The absorbed structure equations take now a simple (sparse) form which
makes the calculations easier.

\begin{figure}[!hbt]
        \begin{center}
        \epsfig{file=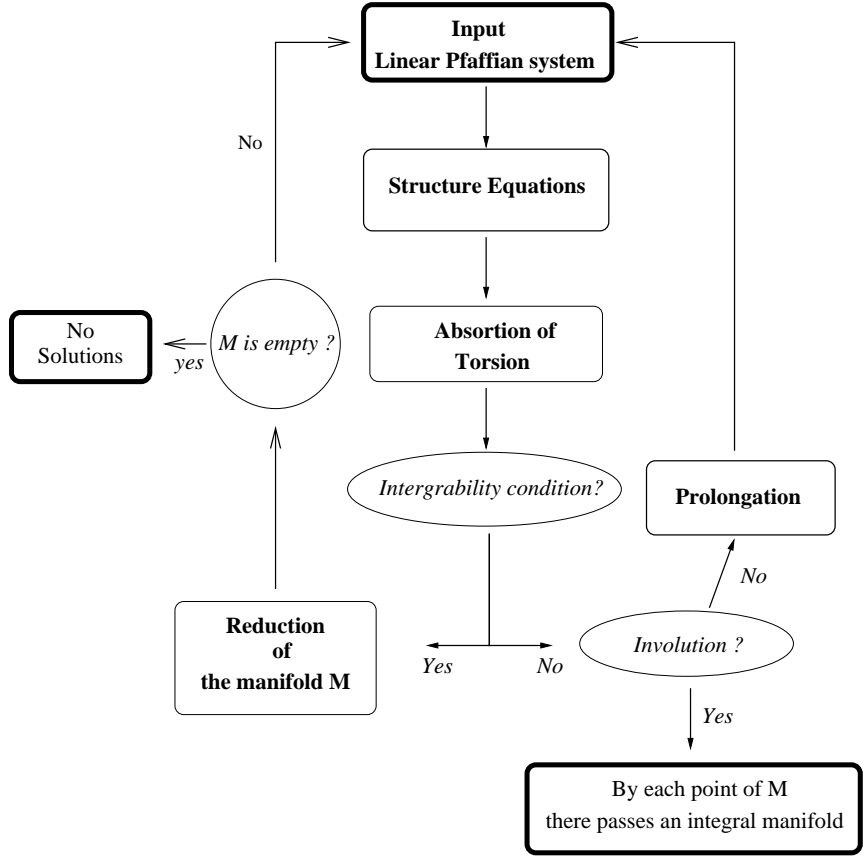, width=11.5cm}
        \end{center}
        \caption{General flowchart}
        \label{algo}
\end{figure}

Before going to the complete  algorithm (see {\sc Fig.}\ref{algo})
 let us sketch two basic concepts which are the
\textit{involution} and the 
\textit{prolongation} (see \cite{Bryant} for rigorous exposition).
By saying that the structure equations are in \textit{involution}  we
simply mean that  by any point of $M$ there passes at least one
ordinary integral  manifold. In practice, involution is checked using  \textit{Cartan characters}.

What about \textit{prolongation}? In the formalism of PDE's systems, a
prolongation consists in  differentiating each equation  w.r.t.
each independent variable.  Thus, one  passes from a system $\E^q
\subset \J^q$ to a  new system $\E^{q+1 } \subset \J^{q+1}$ having the
same solutions.    In the formalism of  linear Pfaffian systems  (and more generally, exterior
differential systems)  the manifold $M$ is replaced by the
\emph{grassmannian}  of the ordinary integral plans. The coordinates
of~$M$  with  the $\lambda^\rho_i$,  which remained arbitrary
after the process of the absorption of torsion, constitute a system of local
coordinates of this grassmannian.  

 Cartan-Kuranishi' theorem \citep{kuranishi} guarantees that the
 algorithm of the figure \ref{algo} stops.  In other words, after a finite number of steps the
 system  is either impossible or in involution
\section{Cartan's method of equivalence}
 \'E. Cartan recasts the problem of local equivalence $\E_f \sim_{\Phi }
 \E_{\bar f}$ into to  the calculation of the integrability conditions of a linear Pfaffian system
\begin{system}\label{Pfaffsys}
 \theta^i_\fbar(\bar\a,\,\bar\x) = \theta^i_f(\a,\,\x),~~(1\leq i \leq m)\\
\theta^{1}_f \wedge \cdots \wedge \theta^{m}_f \neq 0,
\end{system}%
defined on certain manifold $M$ with local coordinates $(\a,\x)$. The change of coordinates $\varphi \in \Phi$ is solution of this
linear Pfaffian system where the  two set of variables
$(\bar\a,\,\bar\x)$ and $(\a,\,\x)$ play a symmetrical role. In this setting,
the integrability conditions  appear under the symmetric form
\begin{equation} \label{invariants}
	 \I_\fbar(\bar\a,\,\bar\x) = \I_f(\a,\,\x)
\end{equation}
The generic function $ \I_f$ is  called \textit{fundamental invariant}.
Every algorithm like {\sc Diffalg} based on the notion of ``ranking'' breaks this symmetry.
Cartan's method, which is an application of the algorithm of the
previous section,  computes the integrability conditions  \pref{invariants} 
by \emph{separately} and \emph{symmetrically} treating
the  $1$-forms $\theta_f$ and $\theta_\fbar$. This divides the number of variables by two.

To fix  the ideas, let us return to the equivalence problem of $2^{nd}$
order ODE introduced in the first section.   
The differential equation $y''=f(x,y,y')$ is equivalent to the Pfaffian system $(\Sigma_f)$
\begin{equation}\
         \d p - f(x,y,p)\ \d x = 0, \quad
        \d y - p\ \d x = 0 \quad \mbox{ with } \quad \d x\wedge \d y \wedge \d p\neq 0.
\end{equation}
In the same way, the equation $\bar y''=\bar f(\bar x,\bar y,\bar y')$
is equivalent  to the Pfaffian system $(\Sigma_{\bar f})$.
Now,  since the first prolongation preserves the module of the contact
  forms,   $(\Sigma_f)$ and $(\Sigma_{\bar f})$ are equivalent
  under  $\varphi \in \Phi$ 
if and only if there exist functions $\a=(a_1, a_2, a_3)$ such that
$$
        \underbrace{
        \begin{pmatrix}
          \d\bar{p}-\bar{f}(\bar x,\bar y,\bar p)\d\bar{x} \\
          \d\bar{y}-\bar{p}\ d\bar{x} \\
          \d\bar{x}
        \end{pmatrix}%
        }_{\omega_\fbar}
        =
\underbrace{
        \begin{pmatrix}
          a_1 & a_2 & 0 \\
          0 & a_3 & 0 \\
          0 & 0 & 1
        \end{pmatrix}%
        }_{ g(\a)}
        \underbrace{
        \begin{pmatrix}
          \d p-f(x,y,p)\ \d x \\
          \d  y-p\ \d x \\
          \d x
        \end{pmatrix}.%
        }_{\omega_f}
$$
The matrices $g(\a)$ form a matrix Lie group, called the \emph{structural} group.
Now, multiplying the two sides by  $g(\a)$ gives 
$$
\underbrace{
        \begin{pmatrix}
          \bar a_1 & \bar a_2 & 0 \\
          0 & \bar a_3 & 0 \\
          0 & 0 & 1
        \end{pmatrix}%
        }_{g(\bar \a)}
        \underbrace{
        \begin{pmatrix}
          \d\bar{p}-\bar{f}(\bar x,\bar y,\bar p)\d\bar{x} \\
          \d\bar{y}-\bar{p}\ d\bar{x} \\
          \d\bar{x}
        \end{pmatrix}%
        }_{\omega_\fbar}
        =
\underbrace{
        \begin{pmatrix}
          a_1 & a_2 & 0 \\
          0 & a_3 & 0 \\
          0 & 0 & 1
        \end{pmatrix}%
        }_{ g(\a)}
        \underbrace{
        \begin{pmatrix}
          \d p-f(x,y,p)\ \d x \\
          \d  y-p\ \d x \\
          \d x
        \end{pmatrix}.%
        }_{\omega}
$$
According to Cartan we define $\theta_f := g(\a)\ \omega_f$ and the
equivalence  problem takes the form \pref{Pfaffsys} where  $\x :=(x,y,p)$.
 We launch Cartan's method as explained in  \citep{Olver1, Gardner,
   Kamran:89, Kamran:90} and we obtain the  involutive structure equations 
\begin{eqnarray}\nonumber
\d\theta^{1} & =&   - \theta^{1} \wedge \theta^{4} + \I_1 \theta^{2} \wedge \theta^{3}, \\ \nonumber
\d\theta^{2} & =&  - \theta^{1} \wedge \theta^{3} - \theta^{2} \wedge \theta^{4}, \\ \nonumber
\d\theta^{3} & =&  0, \\ \nonumber
\d\theta^{4} & =&  \I_2 \theta^{1} \wedge \theta^{2} + \I_3 \theta^{2} \wedge \theta^{3}, \nonumber
\end{eqnarray}
involving the \textit{fundamental} invariants 
\begin{equation*} 
  \I_1(\a,\,\x) = -\Frac{1}{4}(f_{p})^{2} - f_{y} + \Frac{1}{2}D_{x}f_{p}, \quad
  \I_2(\a,\,\x) = \Frac{f_{ppp}}{2 {a_{3}}^{2}}, \quad
  \I_3(\a,\,\x) = \Frac{f_{yp} - D_{x}f_{pp}}{2a_{3}}.
\end{equation*}%
The final invariant 1-forms are
\begin{equation*} \label{ex:thetaforms}
\begin{gathered}
\theta^{1}  =  a_{3}\left( (\Frac{1}{2}f_{p} p-  f ) dx - \Frac{1}{2}f_{p}  dy +  dp \right), \quad 
\theta^{2}  =   a_{3}\left( dy- p dx  \right), \\ 
\theta^{3}  =  dx, \quad 
\theta^{4}  =  ( \Frac{1}{2}f_{p}-\Frac{1}{2}f_{pp} p ) dx + \Frac{1}{2}f_{pp} dy + \Frac{1}{a_{3}}da_{3}.
\end{gathered}
\end{equation*}
 Dual to these forms are  the {invariant derivations}  
\begin{equation*} \label{ex:theta}
\begin{gathered}
X_1 =  \Frac{1}{a_{3}}\Frac{\partial}{\partial p}, \quad
X_2 =  \Frac{1}{a_{3}}\Frac{\partial}{\partial y} + \Frac{1}{2}\Frac{f_{p}}{a_{3}}\Frac{\partial}{\partial p} - \Frac{1}{2}f_{pp}\Frac{\partial}{\partial a_{3}}, \\
X_3 =  \Frac{\partial}{\partial x} + p\Frac{\partial}{\partial y} + f\Frac{\partial}{\partial p} - \Frac{1}{2}f_{p} a_{3}\Frac{\partial}{\partial a_{3}}, \quad
X_4= a_{3}\Frac{\partial}{\partial a_{3}}.
\end{gathered}
\end{equation*}
Thus, the differential of any function $H(x,y,p, a_3)$ can re-expressed as 
$\d H=\sum_{i=1}^4 X_i(H)\ \theta^i.$ 

\section{Differential relations between the fundamental invariants}
The  algebra of invariants associated to a given equivalence problem
is a  {\em differential} algebra generated by the fundamental invariants and closed
under the invariant derivations.
These derivations, which generally do {\emph not} commute, allow us to compute
a \emph{complete} system of invariants from fundamental invariants.
One obtains most of the \emph{syzygies} i.e. the differential relations between the fundamental invariants using
Poincar\'{e} lemma $\d^2 = 0$ where $\d$ denotes
the exterior derivation. This low cost computation does not require
the expression of the invariants in local coordinates or any
elimination which
 is particularly useful since  the invariants can be very big 
(1.1 Mo in the case of ODE's and PDE's systems below).

For the problem (1), the relations provided by  Poincar\'e lemma are
\begin{equation*}
X_1 I_{1} =- I_{3} ,\quad
X_4 I_{1}  =  0,\quad
X_4 I_{2}=- 2 I_{2} ,\quad
X_1 I_{3} =- X_3 I_{2}, \quad
X_4 I_{3} =- I_{3}.
\end{equation*}%
In the particular case when $\bar f=0$,   all invariants vanish and
 thus, according to~\pref{invariants}, the corresponding  invariants of $y''=f(x,y,p)$
must vanish too. Now, since 
 $I_1$ and $I_2$ form a basis of the differential ideal 
generated by the three fundamental invariants, one finds the same conditions (\ref{thm1}).

\subsection{Second order  ODE's systems}
Given the following ODE's system  ($S_F$) :
\begin{system}\label{systutu2}
\ddot{x}^1 = F^1(t,x,\dot{x}), \\
\ddot{x}^2 = F^2(t,x,\dot{x}), \\
\end{system}%
where  $\dot x=(\dot{x}^1, \dot{x}^2)$ denotes the derivative of $x_1$
and $x_2$ according to $t$.
Two systems ($S_F$) and ($S_{\bar F}$) are said to be equivalent 
 under a point transformation if and only if there exist functions $a_1, \dots, a_{15}$ on $\C^3$ in $\C$ such that 
$$ \label{pbtutu1}
        \begin{pmatrix}
        d\dot{\bar{x}}^1 - \overline{F}^1(\bar{t},\bar{x}, \dot{\bar{x}}) d\bar{t} \\
        d\dot{\bar{x}}^2 - \overline{F}^2(\bar{t},\bar{x}, \dot{\bar{x}}) d\bar{t}  \\
        d\bar{x}^1 - \dot{\bar{x}}^1 d\bar{t} \\
        d\bar{x}^2 - \dot{\bar{x}}^2 d\bar{t} \\
        d\bar{t}
        \end{pmatrix}%
        =
        \begin{pmatrix}
        a_1 & a_2 & a_{3} & a_{4} & 0 \\
        a_{5} & a_{6} & a_{7} & a_{8} & 0 \\
        0 & 0 & a_{9} & a_{10} & 0 \\
        0 & 0 & a_{11} & a_{12} & 0 \\
        0 & 0 & a_{13} & a_{14} & a_{15} \\
        \end{pmatrix}
        \begin{pmatrix}
        d\dot{x}^1 - F^1(t,x,\dot{x}) dt \\
        d\dot{x}^2 - F^2(t,x,\dot{x}) dt  \\
        dx^1 - \dot{x}^1 dt \\
        dx^2 - \dot{x}^2 dt \\
        dt
        \end{pmatrix}%
$$
When applied, Cartan's method yields  $88$ fundamental invariants. Without any need to  the explicit expressions of the 88
invariants (over 1 M bytes of memory)  Poincar\'e lemma shows that 
two invariants form a basis the differential
ideal generated by the 88  invariants. If the functions $F^1$ and $F^2$
are identically zero, this two invariants vanish.
\begin{theorem}[\citet{neut:these}]\label{SysEdoThm1}
The system $(S_F)$ is equivalent to the system
$\{ \ddot{\bar{x}}^1 = 0,\ \ddot{\bar{x}}^2 = 0\}$
under a point transformations if and only if
\begin{equation*}
\begin{gathered}
F^2_{\dot{x}^1\dot{x}^1\dot{x}^1}   =  0, \quad
F^1_{\dot{x}^2\dot{x}^2\dot{x}^2}  =  0, \quad
F^2_{\dot{x}^2\dot{x}^2\dot{x}^2} - 3 F^1_{\dot{x}^1\dot{x}^2\dot{x}^2}  =  0, \quad
F^1_{\dot{x}^1\dot{x}^1\dot{x}^1} - 3 F^2_{\dot{x}^1\dot{x}^1\dot{x}^2}  =  0, \\
F^1_{\dot{x}^1\dot{x}^1\dot{x}^2} - F^2_{\dot{x}^1\dot{x}^2\dot{x}^2}  =  0, \quad
2 D_{t}F^1_{\dot{x}^2} - F^1_{\dot{x}^2} F^1_{\dot{x}^1} - F^2_{\dot{x}^2} F^1_{\dot{x}^2} - 4 F^1_{x^2}  =  0, \\
- ({F^2_{\dot{x}^2}})^{2} - 2 D_{t}F^1_{\dot{x}^1} - 4 F^2_{x^2} + 4 F^1_{x^1} + 2\
     D_{t}F^2_{\dot{x}^2} + ({F^1_{\dot{x}^1}})^{2}  =  0, \\
-2 D_{t}F^2_{\dot{x}^1} + F^2_{\dot{x}^2} F^2_{\dot{x}^1} + 4 F^2_{x^1} + F^1_{\dot{x}^1} F^2_{\dot{x}^1}  =  0
\end{gathered}
\end{equation*}
where  
$
 D_t = \frac{\partial}{\partial t } + \dot{x}^1
 \frac{\partial}{\partial x^1 }  + \dot{x}^2 \frac{\partial}{\partial
   x^2 } +  F^1 \frac{\partial}{\partial \dot{x}^1 } + F^2 \frac{\partial}{\partial \dot{x}^2}.
$
\end{theorem}
\subsection{PDE's systems}
Given a system of holomorphic PDE
with two independent variables ${(x^{1 }, x^{2})\in \C^2}$ and one dependent variable $u\in \C$
\begin{equation*}\
   (S_f): \quad \Frac{\partial^2 u}{\partial x^\alpha \partial x^\beta}
        = f_{\alpha \beta} \left( x,\,u,\,\Frac{\partial u}{\partial x} \right) ,\quad
                f_{\alpha \beta}=f_{\beta \alpha} \mbox{ for }\alpha, \beta = 1\ldots 2, \nonumber
\end{equation*}
 If we use the notation
$u' := (u_\alpha)_{1\leq \alpha\leq 2}$ and
$u'' := (u_{\alpha \beta})_{1\leq \alpha \leq \beta \leq 2}$
then  $(S_f)$  reads
$
        u'' = f(x,\, u,\, u').
$  
When $f\equiv 0$ the system is denoted by $(S_0)$. 
 $(S_f)$ and $(S_{\bar f})$
are said to be locally equivalent under a bi-holomorphic transformation  if and only if
\begin{equation}  \label{pb_edp1}
                \begin{pmatrix}
               d\Bar u - \Bar u' d\Bar x\\
               d\Bar x \\
               d\Bar u' - \bar f d\Bar x
        \end{pmatrix}%
        =\begin{pmatrix}
        a & 0 & 0\\
        A & M & 0\\
        B & 0 & N
   \end{pmatrix}
        \begin{pmatrix}
                du - u' dx\\
                dx \\
                du' - f dx
        \end{pmatrix}%
\end{equation}
where $\ a \in \C^*,\ A,\,B \in \C^2,\ M,\,N \in \GL(2, \C)$.
By applying  Cartan's method in S. S. \citet{chern75} way (see also
\citep{fels}),  one obtains  
15  structure equations involving 8 big invariants.
For the  system $(S_0)$, these  invariants vanish. 
\begin{theorem}\label{thm1} The following propositions are equivalent\\
(i) The system $(S_f)$  is equivalent to the  system $(S_0)$ under  bi-holomorphic transformations.\\
(ii) The system $(S_f)$  admits a $15$--dimensional point symmetries  Lie group. \\
(iii) The functions  $f_{\alpha \beta} \mbox{ for }\alpha, \beta = 1\ldots 2$ satisfy
\begin{eqnarray*}\label{cond1}
\Frac{\partial^2 f_{11}}{\partial u_2 \partial u_2} = 0, \quad
\Frac{\partial^2 f_{22}}{\partial u_1 \partial u_1} = 0, \quad
\Frac{\partial^2 f_{12}}{\partial u_2 \partial u_2} -
        \Frac{\partial^2 f_{11}}{\partial u_1 \partial u_2} = 0, \\ \nonumber
\Frac{\partial^2 f_{22}}{\partial u_1 \partial u_2} = 0, \quad
\Frac{\partial^2 f_{11}}{\partial u_1 \partial u_1} -
        4 \Frac{\partial^2 f_{12}}{\partial u_1 \partial u_2}  +
        \Frac{\partial^2 f_{22}}{\partial u_2 \partial u_2} = 0. \nonumber
\end{eqnarray*}
\end{theorem}
\section{Change of coordinates calculation}

One obtains the transformation
$\varphi$ \emph{without} integrating any differential equation when
the symmetry pseudogroup $\S_\fbar \subset \Phi$ of 
the target equation 
$\E_\fbar$ is  zerodimensional.
Indeed, if the function $\varphi$ is the general solution of a
differential system (of non zero order) then
it depends on, at least, one arbitrary constant  and thus (the figure
below)  the symmetry
pseudogroup $\S_\fbar$ is not  zerodimensional.

        \begin{center}
\setlength{\unitlength}{0.00037333in}%
\begingroup\makeatletter\ifx\SetFigFont\undefined%
\gdef\SetFigFont#1#2#3#4#5{%
  \reset@font\fontsize{#1}{#2pt}%
  \fontfamily{#3}\fontseries{#4}\fontshape{#5}%
  \selectfont}%
\fi\endgroup%
{\renewcommand{\dashlinestretch}{30}
\begin{picture}(3311,2233)(0,-10)
\drawline(421,1672)(2349,270)
\drawline(2281.998,297.054)(2349.000,270.000)(2302.612,325.401)
\drawline(2775,1774)(2775,289)
\drawline(2755.845,365.610)(2775.000,289.000)(2794.155,365.610)
\drawline(582,1924)(2491,1924)
\drawline(2419.960,1906.240)(2491.000,1924.000)(2419.960,1941.760)
\put(2700,1849){\makebox(0,0)[lb]{{\SetFigFont{12}{14.4}{\rmdefault}{\mddefault}{\updefault}
\small{$(\xbar_0, \E_\fbar)$}
}}}
\put(-950,1849){\makebox(0,0)[lb]{{\SetFigFont{12}{14.4}{\rmdefault}{\mddefault}{\updefault}
\small{$(\x, \E_f)$}
}}}
\put(1425,2074){\makebox(0,0)[lb]{{\SetFigFont{12}{14.4}{\rmdefault}{\mddefault}{\updefault}
\small{$\varphi_0$}
}}}
\put(2850,874){\makebox(0,0)[lb]{{\SetFigFont{12}{
14.4}{\rmdefault}{\mddefault}{\updefault}
\small{$\sigma\in \S_\fbar$}
}}}
\put(595,874){\makebox(0,0)[lb]{{\SetFigFont{12}{14.4}{\rmdefault}{\mddefault}{\updefault}
\small{$\varphi$}
}}}
\put(2475,-240){\makebox(0,0)[lb]{{\SetFigFont{12}{14.4}{\rmdefault}{\mddefault}{\updefault}
\small{$({\xbar}, \E_\fbar)$}
}}}
\end{picture}
}
        \end{center}


\begin{example}
 Let us go back to the equivalence  with the first Painlev\'{e}
 equation~$(P_I)$ under transformations  $\varphi(x,y)= (\bar x,\bar y)=(x+C,\eta(x,y))$.
We refer the reader to the end of section 4  for the expressions of the fundamental invariants and the invariant derivations.

The specialization of these invariants on the Painlev\'e equation gives
$$ \bar \I_{1}=12\bar y,\quad \bar \I_{2}=\bar \I_{3}=0.$$
According to the equality of the invariants \pref{invariants}, we
deduce that
$\I_2=\I_3=0$
which give the two last equations of \pref{sysp1}. Also, we have

$$
12\bar y=\bar\I_{1}=\I_{1}=-\Frac{1}{4}{f_{p}}^{2} -
f_{y} + \Frac{1}{2} D_{x}f_{p}
$$

and this gives  $\eta(x,y)$, that is the first part of $\varphi$. To
find $C = \bar x -x$, we have
$$
\bar X_3^2 \ \bar \I_1=-72\bar {y}^2-12\bar x
$$
and according again to the equality of the invariants we obtain
\begin{equation*}
C = \bar x - x
  = -\frac{1}{24}\I_1^2-\frac{1}{12}X_3^2 \I_1-x.
\end{equation*}
As byproduct,  the conditions \pref{sysp1} on the function $f$ can be obtained by
expressing that  $C$ is constant i.e. $X_i(C)=0$ for $0\leq i\leq 4$
\begin{eqnarray*}
X_1X_3^2\ \I_1=0, \quad \I_1\ X_2\I_1+X_2X_3^2\ \I_1=0, \quad \I_1\ X_3\I_1+X_3^3\ \I_1-1=0
\end{eqnarray*}
The voluminous formulae \pref{sysp1}  take now a  more compact form, expressed in terms of differential invariants.
\end{example}    

\section{Conclusion}
In this paper we demonstrated  the  superiority of Cartan's
method over direct methods based on differential elimination for
handling equivalence problems. 
Indeed,  we have seen that the 
use of  the  invariant derivations  and the  coding  of the
expressions in terms of invariants  significantly reduce the size of
these expressions. Moreover, in Cartan's method, the frame is
 dynamically adapted (during the computation) using the absorption of
torsion process. This leads to sparse  structure equations and makes 
calculations easier. In addition, this computation is done by
 separately and  symmetrically treating the considered linear Pfaffian
 system. This divides  the number of variables by two.
Also, we have
seen that almost of the syzygies between the fundamental invariants
are obtained using Poincar\'e lemma without any need of the expression of these invariants in
local coordinates (which can take 1 Mo of memory).  
We have gave  new equivalence results,  using our  software which is available at
{\tt www.lifl.fr/\mytilde neut/logiciels}.

\bibliographystyle{elsart-harv}
\bibliography{c}

\end{document}